\documentclass[a4paper,11pt]{article}
\usepackage[T1]{fontenc}
\usepackage[cp1250]{inputenc}
\usepackage{lmodern}

\usepackage[english]{babel}
\usepackage{latexsym}
\usepackage{amsmath}
\usepackage{amsthm}
\usepackage{amssymb}
\usepackage{amsfonts}
\usepackage{stackrel}
\usepackage{dsfont}
\usepackage{tikz}
\usepackage{slashbox}
\usepackage[mathscr]{eucal}
\usepackage{authblk}

\newtheorem{thm}{Theorem}
\newtheorem{lemma}[thm]{Lemma}
\newtheorem{cor}[thm]{Corollary}

\newtheorem{ex}{Example}
\newtheorem{rem}{Remark}
\newtheorem{defin}{Definition}
\newenvironment{sol}{\begin{flushleft} \textbf{Solution}\\* \end{flushleft}}{ \hfill\(\Box\) \\}
\newenvironment{pro}{\begin{flushleft} \textbf{Proof}\\* \end{flushleft}}{\hfill\(\blacksquare\) \\ }

\newcommand{\reals}{\mathbb{R}}
\newcommand{\naturals}{\mathbb{N}}

\newcommand{\complex}{\mathbb{C}}

\newcommand{\eps}{\varepsilon}
\newcommand{\Fred}{\mathscr{F}}
\newcommand{\Hammer}{\mathscr{H}}
\newcommand{\Niem}{\mathscr{N}}
\newcommand{\Urys}{\mathscr{U}}
\newcommand{\Vito}{\mathscr{V}}
\newcommand{\intoper}{\mathscr{P}}

\newcommand{\Addresses}{{
  \bigskip
  \footnotesize

  Krukowski M. (corresponding author), \textsc{Technical Univeristy of \L\'od\'z, \ Institute of Mathematics, \ W\'ol\-cza\'n\-ska 215, \
90-924 \ \L\'od\'z, \ Poland}\par\nopagebreak
  \textit{E-mail address} : \texttt{krukowski.mateusz13@gmail.com}

  \medskip

  Przeradzki B., \textsc{Technical Univeristy of \L\'od\'z,\ Institute of Mathematics, \ W\'ol\-cza\'n\-ska 215, \
90-924 \ \L\'od\'z, \ Poland}\par\nopagebreak
  \textit{E-mail address} : \texttt{przeradz@p.lodz.pl}
}}

\begin{document}
\title{Compactness result and its applications in integral equations}
\author{Mateusz Krukowski$^1$ and Bogdan Przeradzki$^2$}
\affil{$^{1,2}$ Technical University of \L\'od\'z, Institute of Mathematics}
\maketitle

\begin{abstract}
A version of Arzel\`a-Ascoli theorem for $X$ being $\sigma$-locally compact Hausdorff space is proved. The result is used in proving compactness of Fredholm, Hammerstein and Urysohn operators. Two fixed point theorems, for Hammerstein and Urysohn operator, are derived on the basis of Schauder theorem.
\end{abstract}


\section{Introduction}
\noindent

Compactness criteria in typical function spaces not only constitute important results describing properties of these spaces, but they also give a basic tool for investigating the existence of solutions to nonlinear equations of many kinds. The best known criterion is the Arzel\`a-Ascoli theorem that gives necessary and sufficient conditions for compactness in the space of continuous functions defined on a compact space $X$ and taking values in $\reals$ or, more generally, in any finite-dimensional Banach space $E:$ a subset $K\subset C(X,E)$ should consist of equibounded and equicontinuous functions. It is easy to drop the assumption on the dimension of $E$ but then the functions from $K$ should be pointwise relatively compact, i.e. $K_x:=\{ f(x):\, f\in K\}$ is relatively compact for any $x\in X.$ The natural topology in $C(X,E)$ is the topology of uniform convergence given by the norm $\| f\|:=\sup_{x\in X} \|f(x)\|_E.$

If $X$ is not a compact space but only a locally compact one, the Arzela\`a-Ascoli theorem gives a compactness criterion in the space of continuous functions $C(X,Y),$ where $Y$ is a metric space, with the topology of compact convergence, see \cite{Munkres} (or a compact-open topology \cite{Engelking}). A sequence $(f_n)$ tends to $f$ in $C(X,Y),$ if $f_n|_F \to f|_F$ uniformly for each compact subset $F\subset X.$

If one needs the boundedness of this limit $f,$ then he should work in the space of bounded continuous functions $BC(X,Y)$ with its natural topology of uniform convergence. In the case when $Y$ is a Banach space, $BC(X,Y)$ is equipped with the same norm $\| f\| :=\sup \|f(x)\|$ and it is complete. Sufficient conditions for compactness of $K\subset BC(X,Y)$ were known for a long time but they were far from necessary ones and, in fact, they describe subsets of some important subspaces of functions tending to $0$ at infinity (or to any other limits). In \cite{Przeradzki}, one can find sufficient and necessary conditions, but in the case $X=\reals;$ an improved result is in \cite{Stanczy}. 

This paper contains a general case $K\subset BC(X,Y),$ where $X$ is $\sigma$-locally compact and $Y$ is metric and complete (theorem \ref{AAforXY}). This result is applied to linear integral operators (theorem \ref{compactfredholm}) and then to nonlinear integral operators of the Hammerstein (theorem \ref{compacthammerstein}) and Urysohn (theorem \ref{compacturysohn}) type. They give fixed points, solutions of nonlinear integral equations, by using Schauder fixed point theorem. Basic facts about these operators (on compact domains) can be found in \cite{Corduneanu}, \cite{PorterStirling} and \cite{Zemyan}.

\section{Notation}

\noindent
In the paper, we will be mainly preoccupied with the space of bounded and continuous functions. The following definition should clarify any confusion.

\begin{defin}(spaces $C(X,Y)$ and $BC(X,Y)$)\\
Let $(X,\tau_X)$ be a topological space and $(Y,d_Y)$ be a metric space. $C(X,Y)$ is the space of continuous functions with the topology of compact convergence. \\
$BC(X,Y)$ is the space of bounded and continuous functions with the topology of uniform convergence.
\end{defin}

\noindent
Next definition is an extension of the idea found in \cite{Przeradzki}. It is a key concept in theorem \ref{AAforXY}.

\begin{defin}($BC(X,Y)$-extension condition)\\
Let $(X,\tau_X)$ be a topological space and $(Y,d_Y)$ be a metric space. We say that $K \subset BC(X,Y)$ satisfies $BC(X,Y)$-extension condition if for all $\eps > 0$, there exists a compact set $D$ and $\delta > 0$ such that for all $f,g \in K$

$$ d_{BC(D,Y)}(f,g) \leq \delta \ \Longrightarrow \ d_{BC(X,Y)}(f,g) \leq \eps $$

\label{warprzed}
\end{defin}

\noindent
To avoid confusion, we provide the definition of $\sigma$-local compactness, which is used throughout the paper.

\begin{defin}($\sigma$-local compactness)\\
We say that the topological space $(X,\tau_X)$ is $\sigma$-locally compact if it is locally compact and $X = \bigcup_{n=1}^{\infty} S_n$, where $(S_n)_{n=1}^{\infty}$ is a sequence of compact sets. The sequence $(S_n)_{n=1}^{\infty}$ is called the saturating sequence of space $X$.
\end{defin}

\section{Compactness results}

\noindent
The compactness result presented in the paper relies on the following theorem from Munkres 'Topology' \cite{Munkres} (theorem 47.1, page 290) :

\begin{thm}(Arzel\`a-Ascoli for compact Hausdorff space)\\
Let $(X,\tau_X)$ be a locally compact Hausdorff space and $(Y,d_Y)$ be a metric space. Let $K \subset C(X,Y)$, where $C(X,Y)$ is equipped with the topology of compact convergence. $K$ is relatively compact iff it is pointwise relatively compact and equicontinuous at every point of $X$.
\label{AscoliHausdorff}
\end{thm}

\noindent
We use Munkres' theorem on the space $BC(X,Y)$ with topology of uniform (rather than compact) convergence. The improvement is due to the extension condition (definition \ref{warprzed}).

\begin{thm}(Arzel\`a-Ascoli for $\sigma$-locally compact Hausdorff space)\\
Let $(X,\tau_X)$ be $\sigma$-locally compact Hausdorff space and $(Y,d_Y)$ be a complete metric space. The set $K \subset BC(X,Y)$ is relatively compact iff

\begin{description}
	\item[\hspace{0.4cm} (AA1)] $K$ is pointwise relatively compact and equicontinuous at every point of $X$
	\item[\hspace{0.4cm} (AA2)] $K$ satisfies $BC(X,Y)$-extension condition
\end{description}
\label{AAforXY}
\end{thm}
\begin{pro}
'$\Longleftarrow$' \\
Given a sequence from $K$, by theorem \ref{AscoliHausdorff} we choose a subsequence $(f_k)_{k=1}^{\infty}$, which is uniformly convergent on every compact set. For uniform convergence, it suffices to prove that the sequence satisfies Cauchy condition (due to $BC(X,Y)$ being a complete metric space). \\
Let $\eps > 0$. By condition \textbf{(AA2)} we know that there exists a compact set $D$ and $\delta > 0$ such that for all $m,n \in \naturals$ we have

$$ \sup_{x \in D} \ d_Y(f_m(x),f_n(x)) \leq \delta \ \Longrightarrow \ \sup_{x \in X} \ d_Y(f_m(x),f_n(x)) \leq \eps$$

\noindent
Since $(f_k)_{k=1}^{\infty}$ is uniformly convergent on $D$ we obtain that the Cauchy condition is satisfied.\\
'$\Longrightarrow$' \\
\textbf{(AA1)} follows from theorem \ref{AscoliHausdorff}. Assume for the sake of contradiction that \textbf{(AA2)} does not hold. The saturationg sequence $(S_n)_{n=1}^{\infty}$ can be constructed in such a way that $S_n \subset \text{Int}(S_{n+1})$ for every $n \in \naturals$. Hence every compact set $D$ is a subset of $S_n$ for large enough $n$. Negation of \textbf{(AA2)} means that there exists $\eps > 0$ such that for all $n \in \naturals$ there exist $f_n, g_n \in K$ such that

\begin{gather}
\sup_{x \in S_n} \ d_Y(f_n(x),g_n(x)) \leq \frac{1}{n} \ \wedge \ \sup_{x \in X} \ d_Y(f_n(x),g_n(x)) > \eps
\label{negation}
\end{gather}

\noindent
By theorem \ref{AscoliHausdorff}, we choose subsequences $(F_n)_{n=1}^{\infty} \subset (f_n)_{n=1}^{\infty}$ and $(G_n)_{n=1}^{\infty} \subset (g_n)_{n=1}^{\infty}$, which are convergent on compact sets to $F$ and $G$ respectively. The condition (\ref{negation}) implies that for all $m,n \in \naturals$ (without loss of generality $n \geq m$) we have

$$\sup_{x \in S_m} \ d_Y(F_n(x),G_n(x)) \leq \frac{1}{n}$$

\noindent
We conclude that $F \equiv G$, which is a contradiction due to the inequality

$$\sup_{x \in X} \ d_Y(F_n(x),G_n(x)) > \eps $$

\noindent
which holds for every $n \in \naturals$.
\end{pro}

\noindent
We observe that in theorem \ref{AAforXY} the space $X$ can be a subset of finite-dimensional Banach space. In applications it is usually assumed to be closed. The novelty is that we do not assume $X$ to be bounded. We note that if $(Y,d_Y)$ possesses Heine-Borel property then pointwise relative compactness in \textbf{(AA1)} is equivalent to pointwise boundedness.

\section{Application to integral operators}

\subsection{Fredholm operator}

\noindent
In order to avoid misapprehension, we provide the Carath\'eodory conditions which are used throughout the paper. Let $(X,||\cdot||_X), (E,||\cdot||_E)$ be finite-dimensional Banach spaces and $(Y,\Sigma,\mu)$ be a measure space. Let $K : X \times Y \rightarrow B(E)$, where $B(E)$ stands for the space of linear and bounded operators on $E$, satisfy the conditions

\begin{description}
	\item[\hspace{0.4cm} (Car1)] $K(x,\cdot)$ is measurable for every $x \in X$
	\item[\hspace{0.4cm} (Car2)] $K(\cdot,y)$ is continuous for a.e. $y\in Y$
	\item[\hspace{0.4cm} (Car3)] for every $x \in X$ there exists an open neighbourhood $U_x$ and a function $D_x \in L^1(Y)$ such that for every $z \in U_x$ and a.e. $y \in Y$ we have
	
$$||K(z,y)|| \leq D_x(y)$$

\end{description}

\noindent
Observe that \textbf{(Car1)} implies that $||K(x,\cdot)||$ is measurable for every $x \in X$ and \textbf{(Car2)} implies $||K(\cdot,y)||$ is continuous for a.e. $y\in Y$. Hence the function $\Phi : X \rightarrow B(E)$, defined by

\begin{gather*}
\Phi(x) := \int_Y \ ||K(x,y)|| \ d\mu(y)
\label{definitionofphi}
\end{gather*}

\noindent
is continuous for every $x \in X$. Lastly, we adjoin to Carath\'eodory conditions \textbf{(Car1)-(Car3)} the requirement

\begin{description}
	\item[\hspace{0.4cm} (Car4)] the function $K$ satisfies
	
	$$ \sup_{x \in X} \ \int_Y \ ||K(x,y)|| \ d\mu(y) < \infty $$
\end{description}

\noindent
The following lemma is introduced for cosmetic reasons. The main tool used in theorem \ref{compactfredholm} is condition \textbf{(K2)}. However, for simplicity, we provide an equivalent condition \textbf{(K1)}.

\begin{lemma}
Let $(X,||\cdot||_X), (E,||\cdot||_E)$ be finite-dimensional Banach spaces and $(Y,\Sigma,\mu)$ be a measure space. Let $K : X \times Y \rightarrow B(E)$ be a function satisfying conditions \emph{\textbf{(Car1)-(Car4)}}. Then the following conditions are equivalent :

\begin{description}
	\item[\hspace{0.4cm} (K1)] for every $\eps > 0$ and $v \in X, \ ||v|| = 1$ there exist $T_v > 0$ and $L_v \in L^1(Y,B(E))$ such that for every $t \geq T_v$ we have
	
	$$\int_Y \ ||K(tv,y)-L_v(y)|| \ d\mu(y) \leq \eps$$
	
	\noindent
	and moreover $\sup_{||v||=1} T_v < \infty$
	\item[\hspace{0.4cm} (K2)] for every $\eps > 0$ and $v \in X, \ ||v|| = 1$ there exists $T_v > 0$ such that for every $t,s \geq T_v$ we have
	
	$$\int_Y \ ||K(tv,y)-K(sv,y)|| \ d\mu(y) \leq \eps$$
	
	\noindent
	and moreover $\sup_{||v||=1} T_v < \infty$
\end{description}

\label{limitorcauchy}
\end{lemma}
\begin{pro}
\textbf{(K1)} $\Longrightarrow$ \textbf{(K2)} \\
It suffices to note that the triangle inequality

$$||K(tv,y) - K(sv,y)|| \leq ||K(tv,y) - L_v(y)|| + ||L_v(y) - K(sv,y)||$$

\noindent
holds for every $t,s \geq 0, \ v \in X, \ ||v|| = 1$ and a.e. $y \in Y$. \\
\textbf{(K1)} $\Longleftarrow$ \textbf{(K2)} \\
Recognizing \textbf{(K2)} as the Cauchy condition, we obtain the existence of $L_v \in L^1(Y)$ by completeness of the space $L^1(Y, B(E))$.
\end{pro}

\noindent
\begin{rem}
\normalfont The implication \textbf{(K1)} $\Longleftarrow$ \textbf{(K2)} in the proof of lemma \ref{limitorcauchy} does not explicitly construct the limit. In order to obtain a constructive proof, assume additionally

\begin{itemize}
	\item there exists a function $D \in L^1(Y)$ such that for every $x \in X$ and a.e. $y \in Y$, we have
	
	$$||K(x,y)|| \leq D(y)$$
	\item for every $v \in X, \ ||v|| = 1$ and a.e. $y \in Y$ the limit $\lim_{s \rightarrow \infty} \ K(sv, y)$ exists
\end{itemize}

\noindent
Then the limit function in \textbf{(K1)} is given by

$$L_v(y) = \lim_{s \rightarrow \infty} K(sv,y)$$

\noindent
because we can interchange integral sign with the limit in \textbf{(K2)} due to the Lebesgue dominated convergence theorem.
\end{rem}

\begin{thm}
Let $(X,||\cdot||_X), (E,||\cdot||_E)$ be finite-dimensional Banach spaces and $(Y,\Sigma,\mu)$ be a measure space. Let $K : X \times Y \rightarrow B(E)$ satisfy the condtions \emph{\textbf{(Car1)-(Car4)}} and \emph{\textbf{(K1)}}. Then the Fredholm operator $\Fred : L^{\infty}(Y,E) \rightarrow BC(X,E)$ defined by

$$ (\Fred f)(x) := \int_Y \ K(x,y)f(y) \ d\mu(y) $$

\noindent
is linear and compact.
\label{compactfredholm}
\end{thm}
\begin{pro}
We check that the operator is well-defined. The following estimate proves boundedness ($f \in L^{\infty}(Y,E)$ is fixed) of the Fredholm operator

\begin{gather}
\sup_{x \in X} \ ||(\Fred f)(x)|| \leq ||f|| \ \sup_{x \in X} \ \int_Y \ ||K(x,y)|| \ d\mu(y)
\label{boundednessweneed}
\end{gather}

\noindent
For the rest of the proof we fix $\eps > 0$. We prove continuity of the function $(Tf)$ at $x_{\ast} \in X$, where we assume $f \not\equiv 0$ (otherwise continuity is trivial). By conditions \textbf{(Car1)-(Car4)} we know that there exists $\delta > 0 $ such that for every $x \in X, \ ||x - x_{\ast}|| < \delta$ we have

$$ \int_Y \ ||K(x,y)-K(x_{\ast},y)|| \ d\mu(y) \leq \frac{\eps}{||f||}$$

\noindent
Hence, for every $x \in X, \ ||x - x_{\ast}|| < \delta$ we obtain the following inequality

\begin{equation}
	\begin{split}
		||(\Fred f)(x) -  (\Fred f)(&x_{\ast})|| = \bigg|\bigg| \int_Y \ K(x,y)f(y) \ d\mu(y) - \int_Y \ K(x_{\ast}, y)f(y) \ d\mu(y) \bigg|\bigg| \\
								&\leq ||f|| \ \int_Y \ ||K(x,y) - K(x_{\ast},y)|| \ d\mu(y) \leq \eps
	\end{split}
	\label{estimates}
\end{equation}

\noindent
which proves the continuity of $\Fred f$. \\
In order to prove compactness (using theorem \ref{AAforXY}), denote

$$B := \bigg\{ f \in L^{\infty}(Y,E) \ : \ ||f|| < 1 \bigg\}$$

\noindent
$\Fred(B)$ is bounded (in $BC(X,E)$) due to the estimate (\ref{boundednessweneed}). The equicontinuity of $\Fred(B)$ at $x_{\ast}$ follows from the estimate (\ref{estimates}). Hence the condition \textbf{(AA1)} is satisfied. \\
It remains to prove that $\Fred(B)$ satisfies $BC(X,E)$-extension condition, i.e. there exist $T > 0$ and $\delta > 0$ such that for every $f,g \in B$ we have

\begin{gather}
	\sup_{||x|| \leq T} \ \bigg|\bigg| \int_Y \ K(x,y)(f(y)-g(y)) \ dy \bigg|\bigg| \leq \delta \
		\Longrightarrow \ \sup_{x \in X} \ \bigg|\bigg| \int_Y \ K(x,y)(f(y)-g(y)) \ dy \bigg|\bigg| \leq \eps
	\label{whatweneedtoprove}	
\end{gather}

\noindent
By condition \textbf{(K1)} we know that there exists $T>0$ such that for every $v \in X, \ ||v|| = 1$ and every $t,s > T$ we have

$$ \int_Y \ ||K(tv,y)-K(sv,y)|| \ d\mu(y) \leq \frac{\eps}{4} $$

\noindent
For arbitrary $f,g \in B$ denote $h : = f - g$. For every $v \in X, \ ||v|| = 1$ and every $t \geq T$ we have

\begin{equation}
	\begin{split}
		& \bigg| \ \bigg|\bigg| \int_Y \ K(tv,y)h(y) \ d\mu(y) \bigg|\bigg| - \bigg|\bigg| \int_Y \ K(Tv,y)h(y) \ d\mu(y) \bigg|\bigg| \ \bigg| \\
		\leq \bigg|\bigg| \int_Y \ K(tv, y)h&(y)  - K(Tv,y)h(y) \ d\mu(y) \bigg|\bigg| \leq \int_Y \ ||K(tv,y)  - K(Tv,y)|| \ ||h(y)|| \ d\mu(y) \leq \frac{\eps}{2}
	\end{split}
	\label{manymodulus}
\end{equation}

\noindent
By estimate (\ref{manymodulus}), for every $v \in X, \ ||v|| = 1$ and $t \geq T$ we obtain

\begin{gather}
		\bigg|\bigg| \int_Y \ K(tv,y)h(y) \ d\mu(y) \bigg|\bigg| \leq \frac{\eps}{2} + \sup_{ ||x|| \leq T} \ \bigg|\bigg| \int_Y \ K(x,y)h(y) \ d\mu(y) \bigg|\bigg|	
\label{lastestimate}
\end{gather}

\noindent
The estimate (\ref{lastestimate}) gives

$$\sup_{||x|| > T} \ \bigg|\bigg| \int_Y \ K(x,y)h(y) \ d\mu(y) \bigg|\bigg| \leq \frac{\eps}{2} + \sup_{||x|| \leq T} \ \bigg|\bigg| \int_Y \ K(x,y)h(y) \ d\mu(y) \bigg|\bigg|$$

\noindent
Hence it suffices to take $\delta = \frac{\eps}{4}$ for the implication (\ref{whatweneedtoprove}) to be satisfied.
\end{pro}

\begin{cor}
Let $(X,||\cdot||_X)$ be a finite-dimensional Banach space and $(Y,\tau_Y)$ be a topological space with Borel measure $\mu$. Let $K \in BC(X \times Y, \complex)$ satisfy \emph{\textbf{(Car3)-(Car4)}} and \emph{\textbf{(K1)}}. Then the Fredholm operator $\Fred : BC(Y,\complex) \rightarrow BC(X,\complex)$ is linear and compact.
\label{oldcorollary}
\end{cor}

\noindent
In order to present an application of corollary \ref{oldcorollary}, we consider a well-known example of non-compact integral operator $\intoper : BC(\reals,\complex) \rightarrow BC(\reals,\complex)$ defined by

$$(\intoper f)(x) = \int_{\reals} \ e^{-|x - y|} f(y) \ dy$$

\noindent
After small modifications, we are able to obtain a positive result.

\begin{ex}
Let $n \in \naturals$ and $g : \reals^n \rightarrow \reals^n$ be a continuous function such that for every $\eps > 0$ and $v \in X, ||v|| < 1$ there exist $T_v > 0$ and $g_v \in \reals^n$ such that for every $t \geq T_v$ we have $||g(tv) - g_v|| \leq \eps$ and moreover $\sup_{||v||=1} T_v < \infty$. Then the integral operator $\intoper : L^{\infty}(\reals^n,\complex) \rightarrow BC(\reals^n,\complex)$ defined by

$$(\intoper f)(x) = \int_{\reals^n} \ e^{-||g(x) - y||} f(y) \ dy$$

\noindent
is linear and compact.
\label{examplefredholm}
\end{ex}
\begin{sol}
\textbf{(Car3)-(Car4)} follows from the rapid decay of exponential function and the equality

$$\int_{\reals^n} \ e^{||g(x)-y||} \ dy = \int_{\reals^n} \ e^{-||y||} \ dy$$

\noindent
Obviously $e^{-||g_v - \cdot||} \in L^1(\reals^n)$ for every $v \in X, \ ||v|| = 1$. In order to prove \textbf{(K1)} we note the estimate

\begin{gather*}
\bigg| e^{-a} - e^{-b} \bigg| \leq e^{-b} \bigg| e^{-(a-b)} - 1 \bigg| \leq e^{-b} \bigg( e^{|a-b|} - 1 \bigg)
\end{gather*}

\noindent
which holds for $a,b \geq 0$. We set $a = ||g(tv) - y||, \ b = ||g_v - y||$ and integrate over $\reals^n$. For every $v \in X, \ ||v|| = 1$ we obtain

\begin{equation}
	\begin{split}
		\int_{\reals^n} \ \bigg| e^{-||g(tv) - y||} &- e^{-||g_v - y||} \bigg| \ dy \leq \bigg( e^{|| g(tv)-g_v ||} - 1 \bigg) \int_{\reals^n} \ e^{-||g_v - y||} \ dy \\
		&= \bigg( e^{|| g(tv)-g_v ||} - 1 \bigg) \int_{\reals^n} \ e^{-||y||} \ dy
	\end{split}
	\label{exponentestimate}
\end{equation}

\noindent
Fix $\eps >0$ and choose $T>0$ such that for every $v \in X, \ ||v|| =1$ and $t \geq T$ we have

$$ \bigg( e^{|| g(tv)-g_v ||} - 1 \bigg) \int_{\reals^n} \ e^{-||y||} \ dy \leq \eps$$

\noindent
By estimate (\ref{exponentestimate}), we are done.
\end{sol}

\begin{ex}
Let $n \in \naturals$ and $K : \reals^n \times \reals^n \rightarrow \complex$ satisfy \emph{\textbf{(Car1)-(Car4)}} and \emph{\textbf{(K1)}}. Then the Volterra operator $\Vito : L^{\infty}(\reals^n,\complex) \rightarrow BC(\reals^n,\complex)$ defined by

$$ (\Vito f)(x) := \int_{-\infty}^{x_n} \ \ldots \ \int_{-\infty}^{x_1} \ K(x,y)f(y) \ dy_1 \ldots dy_n $$

\noindent
where $x = (x_1, \ldots, x_n)$ and $y = (y_1, \ldots, y_n)$, is linear and compact.
\end{ex}
\begin{sol}
We approximate the Volterra operator $\Vito$ with compact Fredholm operators $(\Fred_m)_{m \in \naturals}$. For every $x,y \in \reals^n$ let

\begin{gather}
K_m(x,y) = K(x,y) \prod_{k=1}^n \left( \mathds{1}_{(-\infty,x_k]}(y_k) + (-my_k + mx_k + 1)\mathds{1}_{(x_k,x_k+\frac{1}{m})}(y_k) \right)
\label{goodkernel}
\end{gather}

\noindent
Kernels $K_m$ satisfy \textbf{(Car1)-(Car4)} and \textbf{(K1)} and hence for every $f \in L^{\infty}(\reals^n,\complex), \ ||f|| < 1$ we have the estimate

\begin{equation}
	\begin{split}
		||\Fred_m f - \Vito f|| \leq \sup_{x \in \reals^n} \
&\int_{x_n}^{x_n+\frac{1}{m}} \ldots \int_{x_1}^{x_1+\frac{1}{m}} \ |K(x,y)| \prod_{k=1}^n \ |-my_k + mx_k + 1| \ dy_1 \ldots dy_n \\
		&\sup_{x \in \reals^n} \ \int_{x_n}^{x_n+\frac{1}{m}} \ldots \int_{x_1}^{x_1+\frac{1}{m}} \ |K(x,y)| \ dy_1 \ldots dy_n
	\end{split}
	\label{approximate1}
\end{equation}

\noindent
The convergence $\Fred_m \rightarrow \Vito$ in the operator norm follows from proposition 3.3.9 (page 112) in \cite{Benedetto}.
\end{sol}

\begin{rem}
\normalfont We could not consider the kernel

\begin{gather*}
x \mapsto K(x,y) \prod_{k=1}^n  \mathds{1}_{(-\infty,x_k]}(y_k)
\label{badkernel}
\end{gather*}

\noindent
in place of (\ref{goodkernel}), beacuse it violates \textbf{(Car2)}.
\end{rem}

\subsection{Hammerstein and Urysohn operator}

\begin{thm}
Let $(X,||\cdot||_X), (E,||\cdot||_E)$ be finite-dimensional Banach spaces, $(Y,\Sigma,\mu)$ be a measure space and $K : X \times Y \rightarrow B(E)$ satisfy the condtions \emph{\textbf{(Car1)-(Car4)}} and \emph{\textbf{(K1)}}. Moreover, let $F : Y \times E \rightarrow E$ satisfy

\begin{description}
	\item[\hspace{0.4cm} (F1)] for every $z \in E$ the function $F(\cdot,z)$ is measurable
	\item[\hspace{0.4cm} (F2)] for every $z \in E$ and $\eps > 0$ there exists $\delta > 0$ such that for a.e. $y \in Y$ and $w \in E$ we have
	
	$$||w - z||\leq \delta \ \Longrightarrow \ ||F(y,w)-F(y,z)||\leq \eps$$
	\item[\hspace{0.4cm} (F3)] there exists a nondecreasing continuous function $\phi : [0,\infty) \rightarrow \reals$ such that for a.e. $y \in Y, \ z\in E$ we have
	
	$$||F(y,z)|| \leq \phi(||z||)$$
\end{description}

\noindent
Then the Hammerstein operator $\Hammer : L^{\infty}(Y,E) \rightarrow BC(X,E)$ defined by

$$ (\Hammer f)(x) := \int_Y \ K(x,y)F(y,f(y)) \ d\mu(y) $$

\noindent
is compact. Moreover, if $X = Y$ and there is $t_{\ast} \in (0,\infty)$ satisfying
	
$$\sup_{x \in X} \ \int_Y \ ||K(x,y)|| \ d\mu(y) \ \phi(t_{\ast}) = t_{\ast}$$

\noindent
then the Hammerstein operator $\Hammer : BC(X,E) \rightarrow BC(X,E)$ has a fixed point.
\label{compacthammerstein}
\end{thm}
\begin{pro}
Observe that the Hammerstein operator $\Hammer$ is the composition of the Fredholm operator $\Fred : L^{\infty}(Y,E) \rightarrow BC(X,E)$ from theorem \ref{compactfredholm} and Nemytskii operator $\Niem : L^{\infty}(Y,E) \rightarrow L^{\infty}(Y,E)$, defined by

$$(\Niem f)(y) := F(y,f(y))$$

\noindent
The fact that Nemytskii operator maps measurable functions to measurable functions follows from \textbf{(F1)-(F2)} and the argument in \cite{Krasnosielskii} (page 349). By \textbf{(F3)}, Nemytskii operator maps bounded sets to bounded sets (in the sense of $L^{\infty}(Y,E)$-norm). Hence, the operator $\Hammer$ maps bounded sets to compact sets. Furthermore, \textbf{(F2)} implies the continuity of $\Niem$ and thus the composition $\Hammer = \Fred \circ \Niem$ is also continuous, which proves that the Hammerstein operator is compact. \\
Assume now, that $X = Y$ and that $\Hammer : BC(X,E) \rightarrow BC(X,E)$. Let $R>0$ and denote

$$B_R := \bigg\{ f \in BC(X,E) \ : \ ||f||\leq R \bigg\}$$

\noindent
For every $f \in B_R$ we have

\begin{equation}
	\begin{split}
		\sup_{x \in X} \ \bigg|\bigg| \int_X \ K(x,y)F(y,&f(y)) \ d\mu(y) \bigg|\bigg| \leq \sup_{x \in X} \ \int_X \ ||K(x,y)||\phi(||f(y)||) \ d\mu(y) \\
		&\leq \sup_{x \in X} \ \int_X \ ||K(x,y)|| \ d\mu(y) \ \phi(R)
	\end{split}
	\label{inequalitywithphi}
\end{equation}

\noindent
\textbf{(F3)} together with (\ref{inequalitywithphi}) implies that $\Hammer(B_{R_{\ast}}) \subset B_{R_{\ast}}$ for some $R_{\ast} >0$. Due to Schauder fixed point theorem, we are done.
\end{pro}

\noindent
We shall show an example of applications of the compactness criterion (Theorem \ref{AAforXY}) to Urysohn integral equations, where the $BC$-extension condition is obtained in another way. Let $\reals_+ := [0,\infty),\ E$ be a finite dimensional Banach space and let $K: \reals_+ \times \reals_+ \times E \to E$ be a continuous map which is uniformly continuous with respect to (w.r.t.) the first variable. Define an integral operator $\Urys : BC(\mathbb{R}_+,E)\to BC(\mathbb{R}_+,E)$ by the formula

$$(\Urys f)(x) := \int_0^{\infty} K(x,y,f(y))\, dy.$$

\noindent
It is well defined if, for any $M>0,$

\begin{equation}
K_M:=\sup_{x\in{\mathbb{R}_+}}\int_0^{\infty} \sup_{|u|\le M} |K(x,y,u)|\, dy<\infty .
\label{finite}
\end{equation}

\noindent
Notice that the uniform continuity of $K$ w.r.t. $x$ is too strong although sufficient condition for $\Urys f$ being continuous but we shall use this assumption in the future. The main assumption (which means that, in a sense, $K$ does not depend on $u$ asymptotically) is the following

\begin{description}
	\item[\hspace{0.4cm} (B)] there exists a function $b\in BC(\reals_+^2,E)$ such that for any $\eps>0$ and $M>0,$ there is $T>0$ with the property
	
$$\int_T^{\infty} |K(x,y,u)-b(x,y)|\, dy <\varepsilon$$

\noindent
for $x\in\mathbb{R}_+,$ $|u|\le M.$
\end{description}

\begin{thm}
Under the above assumptions the Urysohn operator $\Urys$ is compact. If, moreover,

\begin{equation}
\label{limsup}
\limsup_{M\to\infty} \frac{K_M}{M}<1,
\end{equation}

\noindent
then $\Urys$ has a fixed point.
\label{compacturysohn}
\end{thm}
\begin{pro}
First, we prove the equicontinuity of $\Urys f$ for $\| f\|\le M.$ Take $\varepsilon>0$ and $T>0$ from condition \textbf{(B)}. Since $K$ is uniformly continuous w.r.t. $x,$ one can find $\delta>0$ such that

$$|K(x_1,y,u)-K(x_2,y,u)|<\frac{\varepsilon}{T},\quad {\textrm{for}}\quad |x_1-x_2|<\delta,\quad y\in [0,T],\quad |u|\le M.$$

\noindent
Then, for $\|f \| \leq M$ we have

$$|(\Urys f)(x_1)-(\Urys f)(x_2)|\le \int_0^T \ |K(x_1,y,f(y))-K(x_2,y,f(y))|\, dy$$
$$+\int_T^{\infty} \sum_{i=1}^2 \ |K(x_i,y,f(y))-b(x_i,y)|\, dy<3\varepsilon.$$

\noindent
Now, let $f_1,f_2\in BC(\mathbb{R}_+,E)$ be two functions with norms $\le M.$ Taking $T>0$ by condition \textbf{(B)}
and then $\delta>0$ such that, for $x,y\in [0,T],$ $|u_1|,|u_2|\le M$ and $|u_1-u_2|<\delta,$

$$|K(x,y,u_1)-K(x,y,u_2)|<\frac{\varepsilon}{T},$$

\noindent
we get

$$|(\Urys f_1)(x)-(\Urys f_2)(x)|\le \int_0^T \ |K(x,y,f_1(y))-K(x,y,f_2(y))|\, dy$$
$$+\int_T^{\infty} \sum_{i=1}^2 \ |K(x,y,f_i(y))-b(x,y)|\, dy<3\varepsilon.$$

\noindent
provided that $\sup\{|f_1(y) - f_2(y)| \ : \ y \in [0,T]\} \leq \delta$. It follows that $\Urys$ is continuous operator and it satisfies the $BC$-extension condition. Thus, the image of the set of functions with norms bounded by $M$ satisfies the assumptions of Theorem \ref{AAforXY} and $\Urys$ is compact. \\
If condition (\ref{limsup}) holds, then there exists $R>0$ such that $K_R \leq R.$ It follows that the ball $\bar{B}(0,R)$ is mapped into itself by $\Urys$ and a fixed point is given due to the Schauder fixed point theorem.
\end{pro}

\Addresses

\begin{thebibliography}{9}
\bibitem{Benedetto}
	Benedetto J.J., Czaja W. : \textit{Integration and Modern Analysis}, Birkhauser (2009)
\bibitem{Corduneanu}
	Corduneanu C. : \textit{Integral Equations and Applications}, Cambridge
University Press (1991)
\bibitem{Engelking}
  Engelking R. : \textit{General Topology}, PWN Polish Scientific Publishers, Warsaw (1997)
\bibitem{Krasnosielskii}
	Krasnosielskii M.A., Zabreiko P.P., Pustylnik E.I., Sbolevskii P.E. : \textit{Integral operators in spaces of summable functions}, Noordhoff International Publishing, Leyden (1976)
\bibitem{Munkres}
  Munkres J. : \textit{Topology}, Prentice Hall, Inc. (2000)
\bibitem{PorterStirling}
	Porter D., Stirling D. : \textit{Integral equations. A practical treatmment, from spectral theory to applications}, Cambridge University Press (1990)
\bibitem{Przeradzki}
	Przeradzki B. : \textit{The existence of bounded solutions for differential equations in Hilbert spaces}, Annales Polonici Mathematici, LVI.2 (1992)
\bibitem{Stanczy}
	Sta\'nczy R. : \textit{Hammerstein equation with an integral over noncompact domain}, Annales Polonici Mathematici, 69 (1998)
\bibitem{Zemyan}
	Zemyan S.M. : \textit{The Classical Theory of Integral Equations}, Springer Science+Business, LLC (2012)
\end{thebibliography}
\end{document}